\documentclass{amsart}

\newcommand{\ra}{\rightarrow}
\newcommand{\be}{\begin{equation}}
\newcommand{\ee}{\end{equation}}

\newtheorem{theorem}{Theorem}[section]
\newtheorem{lemma}[theorem]{Lemma}
\newtheorem{proposition}[theorem]{Proposition}

\theoremstyle{definition}
\newtheorem{definition}[theorem]{Definition}

\theoremstyle{remark}
\newtheorem{remark}[theorem]{Remark}

\numberwithin{equation}{section}
\begin{document}

\title{Hilbert Space Structure and Positive Operators}

\author{Dimosthenis Drivaliaris}
\address{Department of Financial and Management Engineering\\
University of the Aegean
31, Fostini Str.\\
82100 Chios\\
Greece}
\email{d.drivaliaris@fme.aegean.gr}
\author{Nikos Yannakakis}
\address{Department of Mathematics\\
School of Applied Mathematics and Natural Sciences\\
National Technical University of Athens\\
Iroon Polytexneiou 9\\
15780 Zografou\\
Greece}
\email{nyian@math.ntua.gr}
\commby{}

\begin{abstract}
Let $X$ be a real Banach space. We prove that the existence of an injective, positive, symmetric and
not strictly singular operator from $X$ into its dual implies that
either $X$ admits an equivalent Hilbertian norm or it contains a nontrivially complemented subspace
which is isomorphic to a Hilbert space.
We also treat the non-symmetric case.
\end{abstract}

\maketitle

\section{Introduction}
In this short paper our main aim is to study Banach spaces that contain an isomorphic copy of an infinite dimensional
closed subspace of their dual. This property is directly related to the possibility of a unique representation,
in terms of a bilinear form, of
elements of the dual by elements of the original space. Since this is a common fact
in a Hilbert space it seems  natural to ask if any Hilbert space
structure is present in spaces with this property.
\par
The well known Lax-Milgram theorem
(see \cite{Yosida}) provides sufficient conditions for such a representation in the case of a Hilbert space.
Lin in \cite{Lin} showed
(see theorem \ref{th0} below) that Hilbert space structure is actually necessary, since if the representation
holds then these conditions
imply that the space is isomorphic to a Hilbert space.
\par
In the same spirit, we show that if the operator that carries out the representation is one-to-one, positive and
symmetric then the underlying space is either isomorphic to a Hilbert space or, at least,
it contains a nontrivially complemented subspace which is isomorphic to a Hilbert space. Moreover, we also deal
with the case where the operator is not symmetric by strengthening its positivity
properties.
\par
The history of characterizing Banach spaces that have a Hilbert space structure is a long one and we
do not intend to go through it here. We refer the interested reader to the books of Amir \cite{Amir}
and Istratescu \cite{Istra}.

\section{Preliminaries}
Let $X$ be a real normed space with norm $||\cdot||$, $X^*$ be its dual and
$\langle\cdot,\cdot\rangle$ be their duality product. Let $T:X\ra
X^*$ be a linear operator. We begin with two well known definitions.
\begin{definition}
We say that
$T$ is positive if
$$\langle Tx,x\rangle\geq 0,\text{ for all } x\in X.$$
If in addition $\langle Tx,x\rangle\neq 0$, for $x\neq 0$, we say that $T$ is strictly positive.
\end{definition}
\begin{definition}
We say that $T$ is symmetric if
$$\langle Tx,y\rangle=\langle Ty,x\rangle,\;\text{ for all } x,\;y\in X.$$
\end{definition}
\par
We will use the, already mentioned in the introduction, theorem of Lin \cite{Lin}:
\begin{theorem}
\label{th0}
A Banach space $X$ is isomorphic to a Hilbert space if and only if there exists an isomorphism $T$ from $X$
onto $X^*$ such that
$$m||x||^2\leq \langle Tx,x\rangle\leq M||x||^2,\text{ for all }x\in X$$
where $m$ and $M$ are positive constants.
\end{theorem}
\begin{remark}
It should be noted that the result is still true if $T$ is not assumed to be onto $X^*$.
\end{remark}
\par
For completeness we include the following lemma:
\begin{lemma}
\label{lemma1}
Let $X$ be a real Banach space and $T:X\ra X^*$ be a strictly positive, linear operator. Then:
\begin{itemize}
\item[(i)] $||x||_2=\langle Tx,x\rangle^{\frac{1}{2}}$ is a norm on $X$, there exists a positive constant $c$
such that $||x||_2\leq c||x||$, for all $x\in X$ and
$(X,||\cdot||_2)$ is an inner product space.
\item[(ii)] If $X$ is reflexive, then $T(X)$ is dense in $X^*$.
\end{itemize}
\end{lemma}
\begin{proof}
(i) The strict positivity of $T$ implies that $||x||_2=0\Leftrightarrow x=0$.
The triangle inequality follows from the fact that the positivity of $T$ implies that
$$|\langle Tx,y\rangle+\langle Ty,x\rangle|\leq 2\langle Tx,x\rangle^{\frac{1}{2}}\langle Ty,y\rangle^{\frac{1}{2}},
\text{ for all }x,\;y\in X.$$
The inequality $||x||_2\leq c||x||,\;c>0$  is immediate since every positive operator is bounded. Finally by observing that
$||\cdot||_2$ satisfies the parallelogram law, we get that $(X,||\cdot||_2)$ is an inner product space.\\
(ii) It is a straightforward application of the Hahn-Banach theorem (see Hayden \cite{Hayden}).
\end{proof}

\section{The symmetric case}
In order to prove our main result we need the following Hilbert space characterization:
\begin{proposition}
\label{th1}
A real Banach space $X$ is isomorphic to a Hilbert space if
and only if there exists a  positive and symmetric isomorphism $T:X\ra X^*$.
\end{proposition}
\begin{proof} The necessity is obvious. We show that our claim is
also sufficient. By the positivity and the symmetricity of $T$ we have
that
\be \nonumber |\langle Tx,y\rangle|^2\leq \langle
Tx,x\rangle\langle Ty,y\rangle,\;\;\text{for all}\;\; x,y \in X.
\ee
Let $\varepsilon >0$ and $x\in X$. Then there exists $y\in X$ with $||y||=1$ such that
$$\langle Tx,y\rangle>||Tx||-\varepsilon$$
which implies that
$$\langle Tx,x\rangle^\frac{1}{2}||T||^\frac{1}{2}>||Tx||-\varepsilon.$$
Since $\varepsilon$ is arbitrary we conclude that
$$\langle Tx,x\rangle\geq \frac{1}{||T||}||Tx||^2,\;\;\text{for all } x\in X.$$
Since $T$ is an isomorphism we get
\be
\nonumber
\langle Tx,x\rangle\geq \frac{1}{||T||||T^{-1}||^2}||x||^2,\;\text{ for all }x\in X.
\ee
Therefore there exist $m$ and $M$ positive constants such that
$$m||x||^2\leq \langle Tx,x\rangle\leq M||x||^2,\text{ for all }x\in X.$$
Applying theorem \ref{th0} we conclude that $X$ is isomorphic to a Hilbert space.
\end{proof}
\begin{remark}
If we also assume that $X$ is reflexive, then by using lemma \ref{lemma1}(ii) it can be shown that
$T$ is automatically onto $X^*$. On the other hand, it should be noted that if
$T$ is assumed to be onto $X^*$ then the symmetricity of $T$ alone guarantees that $X$ is reflexive.
For details see Lin \cite{Lin} or Hayden \cite{Hayden}.
\end{remark}
We also need the following
definition:
\begin{definition}
An operator $T:X\ra Y$ is called strictly singular if the restriction of $T$ to any infinite dimensional
subspace of $X$ is not an isomorphism (see \cite{LTz}).
\end{definition}
Our main theorem is the following:
\begin{theorem}
\label{th3} Let $X$ be a real Banach space. If there exists a linear,
injective, symmetric and positive operator $T:X\ra X^*$, which is
not strictly singular, then either $X$ is isomorphic to a Hilbert space
or it contains a nontrivially complemented subspace which is isomorphic to a Hilbert space.
\end{theorem}
\begin{proof} First note that since $T$ is one-to-one, symmetric and positive it is strictly positive.
Hence by lemma \ref{lemma1}(i)
$$||x||_2=\langle Tx,x\rangle^\frac{1}{2}$$
is a norm on $X$ and $(X, ||\cdot||_2)$ is an inner product space. Let
$H$ be the completion of $X$ with respect to $||\cdot||_2$.
Since $T$ is not strictly singular, there exists  an infinite dimensional, closed subspace $M$  of $X$
such that $T(M)$ is closed in $X^*$. By the symmetricity and the positivity of $T$ we have that
$$\langle Tx,x\rangle\geq c||x||^2,\text{ for all }x\in M.$$
Hence, $||\cdot||_2$ is an equivalent norm on $M$  to the initial norm of $X$ and therefore $M$ is also closed in $H$.
Since $H$ is a Hilbert space we have that
\be
\nonumber
H=M\oplus M^\bot
\ee
By noticing that $N=X\cap M^\bot$ is closed in $X$ we get
$$X=M\oplus N.$$
If $M$ has finite codimension then $T$ is an isomorphism on $X$ and by proposition \ref{th1} $X$
is isomorphic to a Hilbert space.
\end{proof}
\begin{remark}
As we have mentioned in the introduction, theorem \ref{th3} implies that a representation of an infinite
dimensional subspace of $X^*$ by a positive,
symmetric and non-degenerate (see \cite{Hayden}) bilinear form is possible
only if either $X$ is itself isomorphic to a Hilbert space or it
contains a nontrivially complemented subspace which is isomorphic to a Hilbert space.
\end{remark}
\begin{remark}
\label{rem}
$l_{1}\oplus l_2$ is an example of a Banach space, not isomorphic to a Hilbert space,
that satisfies the conditions of theorem \ref{th3}.
\end{remark}
\section{The non-symmetric case}
If $T$ is not symmetric, we can still have similar results to those of the previous section
provided that we
strengthen its positivity properties. To this end we need the following definition:
\begin{definition}
We say that a linear operator $T:X\ra X$ is
accretive, if for any $x\in X$, there exists $x^*\in X^*$, such that
$||x^*||^2=||x||^2=\langle x^*,x\rangle$ and $\langle x^*,Tx\rangle\geq 0$.
(see \cite{Pap})
\end{definition}
More suitable for our purposes is the next equivalent definition:
\begin{proposition}
\label{accret}
A linear operator $T:X\ra X$ is
accretive if and only if $$||x+\lambda Tx||\geq ||x||,
\text{ for all } x\in X\text{ and all }\lambda>0.$$
\end{proposition}
(See \cite{Pap})
\begin{remark}
If $X$ is a Hilbert space then, obviously, accretivity is equivalent to positivity.
\end{remark}
In order to proceed we need the following Hilbert space characterization:
\begin{proposition}
\label{th2}
A real Banach space $X$ is isomorphic to a Hilbert space if and only if there exists
a positive isomorphism
$T:X\ra X^*$ onto $X^*$, such that $T^{-1}T^*\tau$ is accretive
(where $\tau$ is the canonical isomorphism of $X$ into $X^{**}$).
\end{proposition}
\begin{proof} The necessity is obvious. To prove the sufficiency of our claim we first note that since $T$ is a
positive isomorphism onto $X^*$, $T+T^*\tau$ is a bounded, linear, symmetric and positive operator from $X$
into its dual. By hypothesis $T^{-1}T^*\tau$ is accretive and hence by proposition \ref{accret},
$I+T^{-1}T^*\tau$ is an isomorphism. But
$$T+T^*\tau=T(I+T^{-1}T^*\tau)$$
and hence $T+T^*\tau$ is also an isomorphism.
Applying proposition \ref{th1} we conclude that $X$ is isomorphic to a
Hilbert space.
\end{proof}
\begin{remark}
(i) If $\;T$ is as in proposition \ref{th1}, then it trivially satisfies the hypotheses of proposition \ref{th2} .\\
(ii) If $X$ is assumed to be reflexive, then the assumption of $T$ being onto $X^*$ can be dropped. \\
(iii) If $X$ is a Hilbert space the requirement that $T^{-1}T^*\tau$ is accretive is equivalent to $T^2$ being
positive.
\end{remark}
Our main theorem for this section reads as follows:
\begin{theorem}
\label{th4}
Let $X$ be a real Banach space. If there exists a
positive isomorphism
$T:X\ra X^*$ onto $X^*$  and an infinite dimensional, closed subspace $M$ of $X$,
such that $T^{-1}T^*\tau|_M$ is accretive, then either $X$ is isomorphic to a Hilbert space or
it contains a nontrivially complemented subspace which is isomorphic to a Hilbert space.
\end{theorem}
\begin{proof} As in the proof of proposition \ref{th2}
$$T+T^*\tau:X\ra X^*$$
is an injective, bounded, linear, symmetric and positive operator
from $X$ into its dual. Arguing as before, we see that
$$(T+T^*\tau)(M)=(T(I+T^{-1}T^*\tau))(M)$$
is closed in $X^*$. Applying theorem \ref{th3} the result follows.
\end{proof}

\end{document}